\documentclass{article}
\usepackage[utf8]{inputenc}
\usepackage[english]{babel}
 
\usepackage[square,sort,comma,numbers]{natbib}
\usepackage{extarrows}
\usepackage{graphicx}
\usepackage{tikz}
\usetikzlibrary{arrows}
\usetikzlibrary{matrix,arrows,decorations.pathmorphing}

\usepackage{graphicx}
\usepackage{tikz}
\usetikzlibrary{arrows}
\usetikzlibrary{matrix,arrows,decorations.pathmorphing}

\usepackage{amsmath}
\usepackage{amssymb}
\usepackage{amsthm}

\usepackage{stmaryrd}
\usepackage[all]{xy}
\usepackage[mathcal]{eucal}
\usepackage{verbatim} %%includes comment environment
\usepackage{fullpage} %%smaller margins
\usepackage{wrapfig}
\usepackage{lipsum}
\usepackage[all]{xy}
\theoremstyle{plain}
\newtheorem{thm}{Theorem}[section]

\newtheorem{lem}[thm]{Lemma}
\newtheorem{prop}[thm]{Proposition}
\newtheorem{cor}[thm]{Corollary}
\newtheorem{claim}[thm]{Claim}
\newtheorem{conj}[thm]{Conjecture}

\theoremstyle{definition}

\theoremstyle{remark}
\newtheorem*{rem}{Remark}

\title{Surjective homomorphisms between surface braid groups}
\author{Lei Chen}

\begin{document}
 \bibliographystyle{alpha}
\maketitle

\begin{abstract}
Let $PB_n(S_{g,p})$ be the pure braid group of a genus $g>1$ surface with $p$ punctures. In this paper we prove that any surjective homomorphism $PB_n(S_{g,p})\to PB_m(S_{g,p})$ factors through one of the forgetful homomorphisms. We then compute the automorphism group of $PB_m(S_{g,p})$, which gives a simpler proof of Irmak--Ivanov--McCarthy \cite[Theorem 1]{ivanov}. Surprisingly, in contrast to the $n=1$ case, any automorphism of $PB_n(S_{g,p})$ for $n>1$ is geometric.

\end{abstract}

\section{Introduction}
Given a surface $S$ and a positive number $n$, we denote the \emph{pure configuration space} of $S$ by
\[
\text{PConf}_n(S):=\{(x_1,...,x_n)\in S^n:x_i\neq x_j\text{ for }i\neq j\}.
\]
There is a natural free action of the permutation group $\Sigma_n$ on PConf$_n(S)$ given by permuting the coordinates and we refer to the corresponding quotient Conf$_n(S):=\text{PConf}_n(S)/\Sigma_n$ as the \emph{configuration space} of $S$. Lastly, denote the \emph{n-strand pure braid group} and the \emph{n-strand braid group} of $S$ by
\[
PB_n(S):=\pi_1(\text{PConf}_n(S))\text{ and }B_n(S):=\pi_1(\text{Conf}_n(S)).
\]
Our goal in this paper is to understand surjective homomorphisms between surface braid groups on different numbers of strands. For example, when $n\ge m>0$, there are natural maps $\text{PConf}_n(S)\to \text{PConf}_m(S)$ forgetting $n-m$ coordinates, which induce \emph{forgetful homomorphisms} $PB_n(S)\to PB_m(S)$. In fact, up to automorphisms, these are the only surjective homomorphisms that arise:
\begin{thm}[surjective homomorphisms $PB_n(S)\to PB_m(S)$]
	Let S be a (possibly noncompact) hyperbolic surface of finite type of genus at least $2$.  For $m,n>0$, every surjective homomorphism $F:PB_n(S)\to PB_m(S)$ factors through some forgetful homomorphism, possibly post-composing with an automorphism of $PB_m(S)$.
	\label{main1}
\end{thm}
In particular, when $m>n$, there is no surjective homomorphism $F:PB_n(S)\to PB_m(S)$. A group $G$ is called Hopfian if any epimorphism of $G$ is an isomorphism. Applying Theorem \ref{main1} in the case $n=m$ gives a proof of the fact that $PB_n(S)$ is Hopfian for $n>0$. Another way to show that $PB_n(S)$ is Hopfian is by showing that $PB_n(S)$ is a finitely generated residually finite group. This comes from an embedding of $PB_n(S)$ inside the automorphism group of a finite generated free group, which is residually finite by Baumslag \cite{MR0146271}. Another consequence of our theorem is the following:
\begin{cor}[surjective homomorphisms $B_n(S)\to PB_m(S)$]
	\label{prop}
Let S be a (possibly noncompact) hyperbolic surface of finite type of genus at least $2$.   For $m>0$ and $n>1$, there is no surjective homomorphism
	\[F:B_n(S)\to PB_m(S).\]
\end{cor}

Historically, the first nontrivial surjective homomorphism between braid groups arose from a classical miracle: ``resolving the quartic". Indeed, let $RQ: \text{Conf}_4(\mathbb{C})\to \text{Conf}_3(\mathbb{C})$ be the map given by
\[RQ(a,b,c,d)=(ab+cd,ac+bd,ad+bc).
\]
By computation, the induced homomorphism on fundamental groups $RQ_*:B_4(\mathbb{C})\to B_3(\mathbb{C})$ is surjective. Theorem \ref{main1} says that there is no such miracle map between pure surface braid groups.

The readers may be wondering why we refer to $RQ_*$ as a miracle. One of the reasons behind this terminology is a result of Lin \cite{lin} saying that there is no surjective homomorphism $B_n(\mathbb{C})\to B_m(\mathbb{C})$ when $n>m$ and $n>4$. To prove this, Lin classified all homomorphisms from $B_n(\mathbb{C})$ to $\Sigma_m$ when $n>m$ extending Artin's \cite{artin} classification of all homomorphisms from $B_n(\mathbb{C})$ to $\Sigma_n$. To get a similar result for surface braid groups, we would need to classify homomorphisms from $B_n(S)$ to $\Sigma_m$, extending Ivanov's \cite[Theorem 1]{ivanovpermutation}  classification of all homomorphisms from $B_n(S)$ to $\Sigma_n$. Based on Theorem \ref{main1}, we have the following conjecture:
\begin{conj}
Let S be a (possibly noncompact) hyperbolic surface of finite type of genus at least two. For $m,n>0$ and $m\neq n$, there is no surjective homomorphism\[
	F:B_n(S)\to B_m(S).\]
\end{conj}

In light of Theorem \ref{main1}, in order to further understand all surjective homomorphisms between surface braid groups, we need to study the automorphism groups of surface braid groups $\text{Aut}(P_n(S))$. To this end, for $g>1$ let $\text{Diff}^{\pm}(S_{g,p,n})$ be the group of diffeomorphisms of $S_g$ fixing two sets of punctures, one with $p$ points and the other with $n$ points including both orientation-preserving and orientation-reversing maps. Let $\text{Mod}^{\pm}(S_{g,p,n}):=\pi_0(\text{Diff}^{\pm}(S_{g,p,n}))$ be the \emph{extended mapping class group} of $S_{g,n,p}$. The following theorem computes $\text{\normalfont Aut}(PB_n(S_{g,p}))$ except when $n=1$ and $p>0$.
\begin{thm}[The automorphism groups of $PB_n(S_{g,p})$ and $B_n(S_{g,p})$]
	\label{main10}
	Assume that $g>1$ and that either $n>1$ and $p>0$ or that $n=1$ and $p=0$. Then
	\[
	\text{\normalfont Mod}^{\pm}(S_{g,p,n})\cong \text{\normalfont Aut}(PB_n(S_{g,p}))\cong \text{\normalfont Aut}(B_n(S_{g,p})).\]

\end{thm}
\begin{rem}
When $n=1$ and $p>0$, the statement of Theorem \ref{main10} is simply false. The group $PB_1(S_{g,p})$ is a free group and there are many isomorphisms of $PB_1(S_{g,p})$ that are not induced from diffeomorphisms. But as Theorem \ref{main10} shows, as long as $n>1$, every automorphism of $PB_n(S_{g,p})$ is induced from a diffeomorphism of $S_{g,p}$.
\end{rem}
It should be mentioned at this point that Theorem \ref{main10} has some predecessors: Irmak--Ivanov--McCarthy \cite[Theorem 1]{ivanov} first computed the automorphism group of $PB_n(S_g)$ and showed that every element is geometric in the sense that it comes from a diffeomorphism of $S_g$. After this work was completed, we also found out that An \cite{MR3488309} obtained Theorem \ref{main10} through a similar method as \cite[Theorem 1]{ivanov}. Moreover, Kida--Yamagata \cite{MR2855817} \cite{MR3125894} showed that every injective homomorphism from a finite index subgroup of $PB_n(S_g)$ to itself is geometric. Nevertheless our method is new and appears to be much simpler than all of the above. In particular, we do not rely on Thurston's theory of surface homeomorphisms and canonical reduction systems. Instead, we use group cohomology and obstruction theory. 
\\
\\
\noindent
{\large\bf Acknowledgements}
The author is grateful to anonymous referees and Maxime Bergeron for suggestions on the writing. She would also like to extend her warmest thanks to Benson Farb for his extensive comments and for his invaluable support from start to finish.
\section{Proof of the classification of surjective homomorphism}
Let $S=S_{g,p}$ be a surface of genus $g$ with $p$ punctures and we denote the \emph{pure configuration space} of $S$ by
 \[
\text{PConf}_m(S):=\{(x_1,...,x_m)\in S^m:x_i\neq x_j\text{ for }i\neq j\}.
\]
In this section, we will use the computations of $H^*(\text{PConf}_m(S_{g,p});\mathbb{Z})$ to prove Theorem \ref{main1}. Most computations here are similar to the computations in Chen \cite[Lemma 3.4]{lei1}. Consider the real codimensional two subspace $\triangle_{ij}$ of $S^m_{g,p}$, defined as 
\[
\triangle_{ij}:=\{(x_1,...,x_m)\in S^m_{g,p}:x_i=x_j\}.
\]
By Poincar\'e duality, the subspace $\triangle_{ij}$ determines a class $[\triangle_{ij}]\in H^2(S_{g,p}^m;\mathbb{Z})$. 

Let $\{a_k',b_k'\}_{k=1}^g$ be a symplectic basis for $H^1(S_{g};\mathbb{Z})$ and $[S_g]'\in H^2(S_{g,p};\mathbb{Z})$ be the fundamental class. Fix a natural embedding $e: S_{g,p}\to S_g$. We define $a_k:=e^*(a_k')$ and $b_k:=e^*(b_k')$ for $k\in \{1,...,g\}$. Let $[S_{g,p}]=e^*([S_g]')\in H^2(S_{g,p};\mathbb{Q})$ be the generator of the fundamental class; when $p>0$, we have that $[S_{g,p}]=0$. Let $p_i:S_{g,p}^m\to S_{g,p}$ be the projection onto the $i$th coordinate. 

Define $H_i:=p_i^*(H^1(S_{g,p};\mathbb{Z}))\subset H^1(S_{g,p}^m;\mathbb{Z})$.
There is a natural embedding $E:\text{PConf}_m(S_{g,p})\to S_{g,p}^m$. By the K\"unneth formula, we have 
\[
H^2(S_{g,p}^m;\mathbb{Z})\cong \bigoplus_{i=1}^{n}\mathbb{Z}p_i^*[S_{g,p}]\oplus \bigoplus_{i\neq j} H_i\otimes H_j.\]
Denote the following composition of maps by $F$: \[
\bigoplus_{i=1}^{m}\mathbb{Z}p_i^*[S_{g,p}]\oplus \bigoplus_{i\neq j} H_i\otimes H_j \xrightarrow[\cong]{\text{K\"unneth}} H^2(S_{g,p}^m;\mathbb{Z})\xrightarrow{E^*} H^2(\text{\normalfont PConf}_m(S_{g,p});\mathbb{Z})
\]
We have the following computations of $H^1(\text{\normalfont PConf}_m(S_{g,p});\mathbb{Z})$ and $H^2(\text{\normalfont PConf}_m(S_{g,p});\mathbb{Z})$.
\begin{lem}\label{comp} Let $g>1$ and $p,m>0$ be integers. \\
(1) We have the following isomorphisms:
\[
H^1(\text{\normalfont PConf}_m(S_{g,p});\mathbb{Z})\xleftarrow[\cong]{E^*}H^1(S_{g,p}^m;\mathbb{Q})\cong \bigoplus_{i=1}^m H_i\]
(2) We have the following exact sequence
\[
0\to \oplus_{1\le i<j\le m}\mathbb{Z}[G_{ij}] \xrightarrow{d} \bigoplus_{i=1}^{n}\mathbb{Z}p_i^*[S_{g,p}]\oplus \bigoplus_{i\neq j} H_i\otimes H_j \xrightarrow{F} H^2(\text{\normalfont PConf}_m(S_{g,p});\mathbb{Z}),
\]
where $d[G_{ij}]=p_i^*[S_{g,p}]+p_j^*[S_{g,p}]+\sum_{k=1}^g (p_i^*a_k \otimes p_j^*b_k-p_i^*b_k\otimes p_j^*a_k)$.
\end{lem}
\begin{proof}
By Totaro \cite[Theorem 1]{totaro}, there is a spectral sequence converging to $H^*(\text{PConf}_m(S_{g,p});\mathbb{Q})$ whose $E_2$ term is a bigraded algebra 
\[
H^{*}(S_{g,p}^m;\mathbb{Z})[G_{ij}],\]
 where $H^r(S_{g,p}^m;\mathbb{Z})$ has degree $(r,0)$ and $G_{ij}$ are generators of degree $(0,1)$ for $1\le i,j\le m$ and $i\neq j$, modulo the following relations:
 \begin{itemize}
 \item $G_{ij}=G_{ji}, (G_{ij})^2=0$,
 \item $ G_{ij}G_{ik}+G_{jk}G_{ji}+G_{ki}G_{kj}=0\text{ for $i,j,k$ distinct}$,
 \item $p_i^*(x)G_{ij}=p_j^*(x)G_{ij}$ for $x\in H^*(S_{g,p};\mathbb{Z})$.  \end{itemize}
 The differential is given by $d_2(G_{ij})=[\triangle_{ij}]$. The following graph is a part of this spectral sequence:
 \begin{center}
\begin{tikzpicture}
\draw [thick] (-0.5,2.5) -- (-.5,0) -- (6,0);
\node at (0.5,0.5) {$\mathbb{Z}$};
\node at (2.5,0.5) {$H^1(S_{g,p}^m;\mathbb{Z})$};
\node (B) at (5,0.5) {$H^2(S_{g,p}^m;\mathbb{Z})$};
\node (A) at (0.5,2) {$\mathbb{Z}[G_{ij}]$};
\node at (2.5,2) {$*$};
\node at (5,2){$*$};
\draw[->] 
  (A) -- (B) node [above,pos=0.5] {$d_2$};
\end{tikzpicture}
\end{center}
By Milnor--Stasheff \cite[Section 11]{CC}, 
\[
[\triangle_{ij}]=p_i^*[S_{g,p}]+p_j^*[S_{g,p}]+\sum_{k=1}^g (p_i^*a_k \otimes p_j^*b_k-p_i^*b_k\otimes p_j^*a_k)\in \bigoplus_{i=1}^{n}\mathbb{Z}p_i^*[S_{g,p}]\oplus \bigoplus_{i\neq j} H_i\otimes H_j \xrightarrow[\cong]{\text{K\"unneth}} H^2(S_{g,p}^m;\mathbb{Z})\]

When $g>1$, elements of the finite set \[
\big\{\sum_{k=1}^g (p_i^*a_k \otimes p_j^*b_k-p_i^*b_k\otimes p_j^*a_k) : 1\le i,j\le m \text{ and } i\neq j \big\}\]
 are linearly independent by the independence of direct sums. Thus $[\triangle_{ij}]$ are linearly independent; i.e. $d_2$ is injective. Then this lemma follows from the convergence of the above spectral sequence.
 \end{proof}
We have the following property about the cup product structure of $H^*(\text{PConf}_m(S_{g,p});\mathbb{Z})$.

\begin{lem}
For two independent elements $x,y \in H^1(\text{\normalfont PConf}_m(S_{g,p});\mathbb{Z})$, if $x\smile y=0$, then there exists $i\in \{1,...,m\}$ such that either $x\in E^*(H_i)$ or $y\in E^*(H_i)$.
\label{crossing}
\end{lem}
\begin{proof}Since $E^*: \bigoplus_{i=1}^m H_i\to H^1(\text{\normalfont PConf}_m(S_{g,p});\mathbb{Z})$ is an isomorphism by Lemma \ref{comp}, we can find $x^i,y^i\in H_i$ for each $1 \le i\le m$ such that $x=E^*(x^1+...+x^m)$ and $y=E^*(y^1+...+y^m)$. The multiplication of $x$ and $y$ satisfies the following:
	\[
	x\smile y=F\big(x^1\smile y^1+...+x^n\smile y^n+\sum_{i\neq j}(x^i\otimes y^j-y^i\otimes x^j)\big)\in \bigoplus_{i=1}^{n}\mathbb{Z}p_i^*[S_{g,p}]\oplus \bigoplus_{i\neq j} H_i\otimes H_j\xrightarrow{F} H^2(\text{PConf}_n(S_{g,p});\mathbb{Z}).
	\]
By $x\smile y=0\in H^2(\text{PConf}_n(S_{g,p});\mathbb{Z})$ and Lemma \ref{comp}, there exists integers $\{k_{i,j}\}$ such that
	\[x^1\smile y^1+...+x^n\smile y^n+\sum_{i\neq j}(x^i\otimes y^j-y^i\otimes x^j)=\sum k_{i,j}d[G_{i,j}] \in 
	\bigoplus_{i=1}^{n}\mathbb{Z}p_i^*[S_{g,p}]\oplus \bigoplus_{i\neq j} H_i\otimes H_j.\] 
	By the independence of all the terms in $\bigoplus_{i=1}^{n}\mathbb{Z}p_i^*[S_{g,p}]\oplus \bigoplus_{i\neq j} H_i\otimes H_j$, we have 
	\begin{equation}
	x^i\otimes y^j-y^i\otimes x^j=k_{i,j}\big(\sum_k (p_i^*a_k \otimes p_j^*b_k-p_i^*b_k\otimes p_j^*a_k)\big)\text{     for all $i,j$}.
	\label{lem22}
	\end{equation}

\begin{claim}
We have that $k_{i,j}=0$ for all $i\neq j$.
\end{claim}
	\begin{proof}
We will prove this claim by contradiction. Assume that $k_{i,j}\neq 0$. In $H_i$, elements $p_i^*a_1,...,p_i^*a_g$ span a subspace $A_i$ and $p_i^*b_1,...,p_i^*b_g$ span a subspace $B_i$. Since $H_i$ is a $\mathbb{Z}$-free module, there is a projection $s_i: H_i\to A_i\oplus B_i$. Let $s_i(x^i)=x_A^i+x_B^i\in A_i\oplus B_i$ and $s_i(y^i)=y_A^i+y_B^i\in A_i\oplus B_i$. By equation \eqref{lem22} and the projection by $s_i\otimes s_j$, we have that 
\[
(x^i_A+x^i_B)\otimes (y^j_A+y^j_B)-(y^i_A+y^i_B)\otimes (x^j_A+x^j_B)=k_{i,j}\big(\sum_{k=1}^g (p_i^*a_k \otimes p_j^*b_k-p_i^*b_k\otimes p_j^*a_k)\big)\in (A_i\oplus B_i)\otimes (A_j\oplus B_j).\]
Since $(A_i\oplus B_i)\otimes (A_j\oplus B_j)=A_i\otimes A_j\oplus A_i\otimes B_j\oplus B_i\otimes A_j\oplus B_i\otimes B_j$, we have 
\begin{itemize}
\item (a) $x^i_A\otimes y_A^j-y_A^i\otimes x^j_A=0$
\item (b) $x^i_B\otimes y_B^j-y_B^i\otimes x^j_B=0$
\item (c) $x^i_A\otimes y_B^j-y_A^i\otimes x^j_B=k_{i,j}(\sum_k p_i^*a_k \otimes p_j^*b_k)$
\item (d) $x^i_B\otimes y_A^j-y_B^i\otimes x^j_A=k_{i,j}(\sum_k p_i^*b_k \otimes p_j^*a_k)$.
\end{itemize}
We claim that $k_{i,j}(\sum_k p_i^*a_k \otimes p_j^*b_k)\in A_i\otimes B_j$ is not a simple tensor of two elements. Assume the contrary that there exists $\lambda_d,\mu_d$ for all $d\in \{1,...,m\}$ such that $z=\sum_d \lambda_d p_i^*a_d$ and $w=\sum_d \mu_d p_i^*b_d$ satisfies that 
 $z\otimes w=k_{i,j}(\sum_k p_i^*a_k \otimes p_j^*b_k)$. Comparing the coefficient of $p_i^*a_k\otimes p_j^*b_k$, we know that $\lambda_k\neq 0$ and $\mu_k\neq 0$ for any $k$. Since $g>1$, the coefficient of $p_i^*a_1\otimes p_j^*b_2$ is nonzero, which is a contradiction.
 
By equation (c) and the fact that $\sum_k p_i^*a_k \otimes p_j^*b_k\in A_i\otimes B_j$ is not a simple tensor, we know that $x_A^i,y_B^j,y_A^i$ and $x_B^j$ are all nonzero. For the same reason, $x_B^i,y_A^j,y_B^i$ and $x_A^j$ are also all nonzero. Equation (a) says that there exists $\mu$ such that $x_A^i=\mu y_A^i$ and $\mu y_A^j=x_A^j$; equation (b) says that there exists $\mu'$ such that  $x_B^i=\mu' y_B^i$ and $\mu' y_B^j=x_B^j$. Therefore 
\[
x^i_A\otimes y_B^j-y_A^i\otimes x^j_B=\mu y_A^i\otimes y_B^j-y_A^i\otimes \mu'y_B^j=(\mu-\mu')y_A^i\otimes y_B^j,\]
which contradicts the fact that $\sum_k p_i^*a_k \otimes p_j^*b_k\in A_i\otimes B_j$ is not a simple tensor. Thus the claim holds.
	\end{proof}

	 Therefore $x^i\otimes y^j-y^i\otimes x^j=0\in H_i\otimes H_j$. Assume that $x\notin H_i$ for any $i$; i.e. two coordinates of $x$ are nonzero. Without loss of generality, we assume that $x^1\neq 0$ and $x^2\neq 0$. We break the proof into the following cases.
	\begin{itemize}
\item Case 1) $y^1\neq 0$ and $y^1$ is not proportional to $x^1$. Then $x^1\otimes y^j=y^1\otimes x^j\in H_1\otimes H_j$ implies that $y^j=0$ and $x^j=0$ for all $j\neq 1$, which contradicts to the assumption that $x_2\neq 0$.
	\item Case 2) $y^1\neq 0$ and $x^1=\mu y^1$ for $\mu\neq 0$. Then $x^1\otimes y^j=y^1\otimes x^j\in H_1\otimes H_j$ implies that $x^j=\mu y^j$ for all $j$. Therefore
	\[
	x=x^1+...+x^m=\mu y^1+...+\mu y^m=\mu y,\]
	which contradicts the fact that $x$ and $y$ are independent.
	\item Case 3) $y^1=0$. Then $x^1\otimes y^j=y^1\otimes x^j\in H_1\otimes H_j$ implies that $y^j=0$ for all $j$, which contradicts to the fact that $x$ and $y$ are independent.\qedhere 	\end{itemize}\end{proof}
We define $PB_n(S_{g,p}):=\pi_1(\text{PConf}_n(S_{g,p}))$ and let $P_i:=(p_i\circ E)_*: PB_n(S_{g,p})\to \pi_1(S_{g,p})$ be the induced map on the fundamental groups of the projection to the $i$th coordinate. 
\begin{prop}\label{fg1}
For $g>0$ and $p,n\ge 0$, the space
$\text{PConf}_n(S_{g,p})$ is a $K(\pi,1)$-space and Ker$(P_i)$ is finitely generated.
\end{prop}
\begin{proof}
We will prove that $\text{PConf}_n(S_{g,p})$ is a $K(\pi,1)$-space for any $p\ge 0$ by induction on $n$. For $n=1$, the space $\text{PConf}_{n}(S_{g,p})=S_{g,p}$ is a $K(\pi,1)$-space for any $p$. Assume that $\text{PConf}_{n-1}(S_{g,p})$ is a $K(\pi,1)$-space for any $p$ for $n>1$. The forgetful map $p_i\circ E$ gives the following fiber bundle
\begin{equation}
\text{PConf}_{n-1}(S_{g,p+1})\to \text{PConf}_n(S_{g,p})\xrightarrow{p_i\circ E} S_{g,p}.
\label{fb}
\end{equation}
Since both $\text{PConf}_{n-1}(S_{g,p+1})$ and $S_{g,p}$ are $K(\pi,1)$-spaces, we have that $\text{PConf}_n(S_{g,p})$ is a $K(\pi,1)$-space. Then the induction axiom implies that $\text{PConf}_n(S_{g,p})$ is a $K(\pi,1)$-space for any $n$. Fiber bundle \eqref{fb} induces the following short exact sequence on fundamental groups:
\begin{equation}
1\to \text{PB}_{n-1}(S_{g,p+1})\to \text{PB}_n(S_{g,p})\xrightarrow{P_{i}} \pi_1(S_{g,p})\to 1.
\label{fb1}
\end{equation}
By the short exact sequence \eqref{fb1}, we obtain Ker$(P_{i})=PB_{n-1}(S_{g,p+1})$. Using \eqref{fb1}, it is straightforward to prove that $PB_n(S_{g,p})$ is finitely generated by induction on $n$. Therefore Ker$(P_{i})$ is finitely generated.
\end{proof}
Now, we proceed to the key lemma.
\begin{lem}For $g>1$ and $p,n\ge 0$, a homomorphism
\[
R: PB_n(S_{g,p})\to \pi_1(S_{g,p})
\] 
either factors through $P_i$ for some $i\in \{1,...,n\}$ or has cyclic image.
\label{cyclic}
\end{lem}

\begin{proof}
The proof of this lemma uses the same idea as F.E.A. Johnson \cite{FEAJohnson}. The method can also be found in Salter \cite[Lemma 3.3 and 3.4]{NickFibering}. We use group cohomology in what follows. Since all spaces we consider are $K(\pi,1)$-spaces, we freely use the computation of the cohomology of the spaces as the cohomology of the corresponding groups. By the classification of subgroups of $\pi_1(S_{g,p})$, if Image$(R)\ncong \mathbb{Z}$, then Image$(R)$ is either a free group $F_k$ with $k>1$ or a surface group $\pi_1(S_h)$ such that $h\ge g>1$. In both cases, there are independent elements $x,y\in H^1(\text{Image}(R);\mathbb{Z})$ such that $x\smile y=0$. Denote by $S:PB_n(S_{g,p})\to \text{Image}(R)$ the map to the image of $R$, which is surjective by definition. Then $S^*(x),S^*(y)\in H^1(PB_n(S_{g,p});\mathbb{Q})$ are independent and $S^*(x)\smile S^*(y)=0$. By Lemma \ref{crossing}, we have either $S^*(x)\in P_i^*\big(H^1(\pi_1(S_{g,p});\mathbb{Z})\big)$ or $S^*(y)\in P_i^*\big(H^1(\pi_1(S_{g,p});\mathbb{Z})\big)$ for some $i$. Without loss of generality, assume that $S^*(x)=P_i^*(x')$. We have the following commutative diagram by the identification $H^1(\_\_;\mathbb{Z})\cong \text{Hom}(\_\_,\mathbb{Z})$,
\[
\xymatrix{
	PB_n(S_{g,p})\ar[r]^-S\ar[d]^{P_i} & \text{Image}(R)\ar[d]^x\\
	\pi_1(S_{g,p})\ar[r]^{x'} & \mathbb{Z}.}
\]
By Proposition \ref{fg1}, the group $\text{Ker}(P_{i})$ is a finitely generated normal subgroup of $PB_n(S_{g,p})$. Since $S$ is surjective, the image $S\big(\text{Ker}(P_{i})\big)$ is also a finite generated normal subgroup of Image$(R)$. However every finitely generated normal subgroup of Image$(R)$ is either of finite index or trivial; see e.g. F.E.A. Johnson \cite[Property (D6)]{FEAJohnson}. If $S\big(\text{Ker}(P_{i})\big)<\text{Image}(R)$ is of finite index, then after composing with $x$, the image $x\circ S\big(\text{Ker}(P_{i})\big)$ will be of finite index in $\mathbb{Z}$. This is a contradiction because $x\circ S\big(\text{Ker}(P_{i})\big)=x'\circ P_i\big(\text{Ker}(P_{i})\big)=\{1\}$. Therefore $S\big(\text{Ker}(P_{i})\big)=1$; i.e. $S$ factors through $P_i$.
\end{proof}
Now we are ready to prove Theorem \ref{main1} saying that for $m,n>0$ and $g>1$, every surjective homomorphism $F:PB_n(S_{g,p})\to PB_m(S_{g,p})$ factors through some forgetful homomorphism, possibly post-composing with an automorphism of $PB_m(S_{g,p})$.

\begin{proof}[Proof of Theorem \ref{main1}]
We will prove Theorem \ref{main1} by induction on $m$. There are two things to be proved: every surjective homomorphism $F: PB_n(S_{g,p})\to PB_m(S_{g,p})$ factors through a forgetful homomorphism for any $m$ and $n$ and that $PB_m(S_g)$ is Hopfian for any $m$. The case $m=1$ is a result of Lemma \ref{cyclic} and the classical fact that $\pi_1(S)$ is Hopfian; see e.g. Hempel \cite{MR0295352}. We assume that when $m<k$, Theorem \ref{main1} is true. For $m=k$, let $f: PB_n(S_{g,p})\to PB_k(S_{g,p})$ be a surjective homomorphism. By post-composing with a projection $pr:PB_k(S_{g,p})\to PB_{k-1}(S_{g,p})$, we obtain a new surjective homomorphism $pr\circ f:PB_n(S_{g,p})\to PB_{k-1}(S_{g,p})$. By the inductive hypothesis, $pr\circ f$ factors through some forgetful map, possibly post-composing with an automorphism of $PB_{k-1}(S_{g,p})$. Therefore we have the following commutative diagram:
\begin{equation}\label{bigdiagram}
\xymatrix{
1\ar[r] &PB_{n-k+1}(S_{g,p+k-1})\ar[r]\ar[d]^r & PB_n(S_{g,p})\ar[r]^{\text{forget}}\ar[d]^f & PB_{k-1}(S_{g,p})\ar[r]\ar[d]^p & 1\\
1\ar[r] &PB_1(S_{g,p+k-1})\ar[r] & PB_k(S_{g,p})\ar[r]^{pr} & PB_{k-1}(S_{g,p})\ar[r] & 1}.
\end{equation}
Since $PB_{n-k+1}(S_{g,p+k-1})\to PB_1(S_{g,p+k-1})$ factors through a forgetful homomorphism, which is the $m=1$ case of the theorem, we show that $f$ factors through a forgetful homomorphism. When $n=k$, both $r$ and $p$ are isomorphisms because of the inductive hypothesis. Therefore $f$ is an isomorphism by the five lemma. We finish the proof by the induction axiom.
\end{proof}
Define $B_n(S):=\pi_1(\text{PConf}_n(S)/\Sigma_n)$, where $\Sigma_n$ acts on PConf$_n(S)$ by permuting coordinates. Now we are ready to prove Corollary \ref{prop} saying that there is no surjective homomorphism from $B_n(S_{g,p})$ to $PB_m(S_{g,p})$ when $g>1, n>1$ and $m>0$. 
\begin{proof}[Proof of Corollary \ref{prop}]
Set $S=S_{g,p}$. We plan to prove that there is no surjective homomorphism $f:B_n(S_{g,p})\to PB_m(S_{g,p})$ by contradiction. Assume the opposite that there is a surjective homomorphism  $f:B_n(S_{g,p})\to PB_m(S_{g,p})$. By projecting to some coordinate, we obtain a  surjective homomorphism $F:B_n(S)\to \pi_1(S)$. By post-composing with the embedding $p_*: PB_n(S)\to B_n(S)$, which is induced by the projection $p: \text{PConf}_n(S)\to \text{Conf}_n(S)$, we obtain a map $R=F\circ p_*: PB_n(S)\to \pi_1(S)$. Since $F$ is surjective, Im$(R)$ is an index $n!$ subgroup $H<\pi_1(S)$, where $H=\pi_1(S')$ is the fundamental group of an $n!$ cover $S'$ of $S$. We have the following inequalities of the first Betti numbers of $S$ and $S'$:
\begin{align*}
b_1(S')&\ge 1-\chi(S')    && \text{$b_1(S')=b_0(S')+b_2(S')-\chi(S')\ge 1-\chi(S')$} \\
&=1-n!\chi(S)      && \text{$S'$ is $n!$-cover of $S$} \\
&\ge 1-2\chi(S) &&  \text{$\chi(S)$ is negative}\\
&\ge 3-\chi(S)    && \text{$g>1$ implies that $\chi(S)\le -2$}\\
&\ge 1+b_1(S)&& \text{$b_1(S)=b_0(S)+b_2(S)-\chi(S)\le 2-\chi(S)$.}
\end{align*}
Since $R$ does not have a cyclic image, we know that $R=Q\circ P_i$ for some $i\in \{1,...,n\}$ and some map $Q:\pi_1(S)\to \pi_1(S)$ by Lemma \ref{cyclic}. Since $P_i$ is surjective, we know that Im$(Q)=\text{Im}(R)=H$. However $Q$ cannot be surjective onto the image $\pi_1(S')$ because $b_1(S)<b_1(S')$, which is a contradiction.
\end{proof}

\section{The automorphism group of PB$_n(S_{g,p})$}

In this section, we will compute the automorphism group of PB$_n(S_{g,p})$. The key point is to use the existence of pseudo-Anosov elements in the point-pushing subgroup. Before the proof of the result, we will introduce $3$ classical results we will use in the proof. Firstly, we have the following result of Handel--Thurston \cite[Lemma 2.2]{handel}.
\begin{thm}[\cite{handel}]
A pseudo-Anosov element of the mapping class group does not fix any nonperipheral isotopy class of curves (including nonsimple curves).
\label{thur}
\end{thm}
Another ingredient is Kra's construction \cite{Kra}. Let $S$ be a finite type surface possibly with punctures. The \emph{extended mapping class group} Mod$^\pm(S)$ is defined to be the group of isotopy classes of diffeomorphisms of $S$ fixing the punctures as a set. Later, we will define other types of extended mapping class groups by specifying exactly how they preserve the punctures. Let $b\in S$. Denote by Mod$^\pm(S,b)$ the extended mapping class group of $S$ fixing a point $b$; in particular Mod$^\pm(S,b)$ fixes the punctures of $S$ as a set and the point $b$. The following is the Birman exact sequence for $S$ (see e.g. Farb--Margalit \cite[Section 4.2]{BensonMargalit}):
\[
1\to \pi_1(S,b)\xrightarrow{Push}\text{Mod}^\pm(S,b)\to \text{Mod}^\pm(S)\to 1.
\]
We say that a nontrivial element $\gamma\in \pi_1(S,b)$ \emph{fills} $S$ if the curve representing $\gamma$ intersects every essential simple closed curve on $S$.
\begin{thm}[Kra's construction \cite{Kra}]
Let $S$ be a finite type surface possibly with punctures. Let $\gamma\in \pi_1(S,b)$. The mapping class $Push(\gamma)\in \text{\normalfont Mod}^\pm(S,b)$ is pseudo-Anosov if and only if $\gamma$ fills $S$.
\label{kra}
\end{thm}
The third ingredient is the following punctured Dehn--Nielsen--Baer theorem; e.g. see e.g. Farb--Margalit \cite[Theorem 8.8]{BensonMargalit}. For a group $G$, denote by Out$(G)$ the \emph{outer automorphism group} of $G$; i.e. Out$(G)=\text{Aut}(G)/\text{Inn}(G)$, where $\text{Aut}(G)$ denotes the automorphism group of $G$ and Inn$(G)$ denotes the group subgroup of $\text{Aut}(G)$ consisting of conjugations.

\begin{thm}[Punctured Dehn--Nielsen--Baer Theorem]\label{DNB}
Let $S$ be a finite type surface possibly with punctures. Let {\normalfont Out}$^*(\pi_1(S))$ (resp. {\normalfont Aut}$^*(\pi_1(S))$) be the subgroup of {\normalfont Out}$(\pi_1(S))$ (resp. {\normalfont Aut}$(\pi_1(S))$) consisting of elements that leave invariant the set of conjugacy classes in $\pi_1(S)$ of simple closed curves surrounding individual punctures. Then the natural maps
	\[
	\text{\normalfont Mod}^{\pm}(S)\to \text{\normalfont Out}^*(\pi_1(S)) \text{ and  } 	\text{\normalfont Mod}^{\pm}(S,b)\to \text{\normalfont Aut}^*(\pi_1(S))
	\]
	are isomorphisms.	
\end{thm}
Given $f\in \text{Mod}^{\pm}(S_{g,p+n})$, there is an induced action $f_*$ on the fundamental group $\pi_1(S_{g,p+n})$ up to conjugation. Therefore we have an injective homomorphism $\rho:\text{Mod}^{\pm}(S_{g,p+n})\to \text{Out}(\pi_1(S_{g,p+n}))$ given by $\rho(f)=f_*$. By post-composing with the embedding of the braid point-pushing subgroup $PB_n(S_{g,p})< \text{Mod}^{\pm}(S_{g,p+n})$, we obtain an embedding
\[
\theta: PB_n(S_{g,p})< \text{Mod}^{\pm}(S_{g,p+n})\to \text{Out}(\pi_1(S_{g,p+n})).
\]
Let Mod$^\pm(S_{g,p,n})$ be the extended mapping class group that fixes two sets of punctures, one with $p$ points and the other with $n$ points. The following proposition computes the normalizer of $PB_n(S_{g,p})<\text{Out}(\pi_1(S_{g,p+n}))$, which is the main ingredient in proving Theorem \ref{main10}. 

\begin{prop}
For $n>0$, $g>1$ and $p\ge 0$, the normalizer of $PB_n(S_{g,p})<\text{Out}(\pi_1(S_{g,p+n}))$ is $\text{\normalfont Mod}^{\pm}(S_{g,p,n})$.
\end{prop}
\begin{proof}
Let $R\in \text{Out}(\pi_1(S_{g,n+p}))$ be an element in the normalizer of $PB_n(S_{g,p})$. Then $R$ acts as an automorphism $A$ on $PB_n(S_{g,p})$ by conjugation; i.e. $R\theta(e)R^{-1}=\theta(A(e))$ for $e\in PB_n(S_{g,p})$. This gives the following equation
\begin{equation}
R\theta(e)=\theta(A(e))R \in \text{Out}(\pi_1(S_{g,p+n})).
\label{prop34}
\end{equation}
We claim that for any curve $c$  surrounding a puncture, $R(c)$ is also a curve surrounding a puncture. This implies that $R\in \text{\normalfont Mod}^{\pm}(S_{g,p+n})$ by Theorem \ref{DNB}.

Let $c$ be a simple closed curve surrounding a puncture. Since $PB_n(S_{g,p})$ fixes all punctures, we have $\theta(e)(c)=c$ and $\theta\big(A(e)\big)(c)=c$. By equation \eqref{prop34}, we have that $\theta\big(A(e)\big)(R(c))=R\big(\theta(e)(c)\big)=R(c)$ for any $e\in PB_n(S_{g,p})$. However, because of Theorem \ref{kra}, we know that that there is a pseudo-Anosov element $e'$ in $PB_n(S_{g,p})$. Set $e=A^{-1}(e')$. Then $\theta\big(A(e))(R(c)\big)=\theta(e')(R(c))$. Therefore $\theta(e')\big(R(c)\big)=R(c)$, which implies that $R(c)$ is peripheral by Theorem \ref{thur}.

What remains to be proven is that $R\in \text{Mod}^{\pm}(S_{g,p,n})$. By the \emph{generalized Birman exact sequence} (see e.g. Farb--Margalit \cite[Theorem 9.1]{BensonMargalit})
\[
1\to PB_{n+p}(S_g)\to \text{Mod}^{\pm}(S_{g,p+n})\to \text{Mod}^{\pm}(S_{g})\times \Sigma_n\to 1,
\]
$R$ induces an automorphism $A'$ of $PB_{n+p}(S_g)$ by conjugation. We have the following exact sequence given by forgetting the $n$ punctures 
\begin{equation}
1\to PB_n(S_{g,p})\to PB_{n+p}(S_g)\to PB_p(S_g)\to 1.
\label{forget}
\end{equation}
Since $A'$ preserves the subgroup $PB_n(S_{g,p})\unlhd PB_{n+p}(S_g)$, it should induce an action on each term of \eqref{forget}. As a result, $R$ fixes the $p$ punctures as a set, which implies that $R$ should also fix the other $n$ punctures as well. Thus $R\in \text{Mod}^{\pm}(S_{g,p,n})$.
\end{proof}
Let Mod$^\pm(S_{g,p,n-1,1})$ be the extended mapping class group that fixes three sets of punctures, one with $p$ points, one with $n-1$ points and the last one with one point. Let $\phi: \text{Mod}^{\pm}(S_{g,p+n-1,1})\to \text{Aut}\big(\pi_1(S_{g,p+n-1})\big)$ be the natural embedding induced by conjugation on the normal subgroup $\pi_1(S_{g,p+n-1})\unlhd \text{Mod}^{\pm}(S_{g,p+n-1,1})$. By a similar argument, we obtain the following:
\begin{prop}	\label{haha} For $g>1,n>1,p\ge 0$, the normalizer of

\[
PB_n(S_{g,p})<\text{\normalfont Mod}^{\pm}(S_{g,p+n-1,1})\xrightarrow{\phi}\text{\normalfont Aut}\big(\pi_1(S_{g,p+n-1})\big)\text{ is $\text{\normalfont Mod}^{\pm}(S_{g,p,n-1,1})$}.\]
\end{prop}

We are now ready to prove Theorem \ref{main10}; that is $ \text{Mod}^\pm(S_{g,p,n})\cong \text{Aut}(PB_n(S_{g,p}))$.

\begin{proof}[\bf Proof of Theorem \ref{main10}]
First of all, from another version of the \emph{generalized Birman exact sequence} (see e.g. Farb--Margalit \cite[Theorem 9.1]{BensonMargalit})
\[
1\to PB_n(S_{g,p})\to \text{Mod}^\pm(S_{g,p,n})\to \text{Mod}^{\pm}(S_{g,p})\times \Sigma_n\to 1, \]
there is a map 
\[
C: \text{Mod}^\pm(S_{g,p,n})\to \text{Aut}(PB_n(S_{g,p}))
\]
given by conjugating the subgroup $PB_n(S_{g,p})$. Let $(b_1,...,b_n)\in \text{PConf}_n(S_{g,p})$ be a base point. The mapping class group $\text{Mod}^\pm(S_{g,n,p})$ fixes $\{b_1,...,b_n\}$ as a set. Let $T:\text{Mod}^{\pm}(S_{g,p,n})\to \Sigma_n$ be the homomorphism recording to the permutation of the $n$ points $\{b_1,...,b_n\}$. Let $f\in \text{Mod}^\pm(S_{g,n,p})$. A geometric representation $F\in \text{Diff}^\pm(S_{g,p,n})$ of $f$ induces a map $f_n: \text{PConf}_n(S_{g,p})\to \text{PConf}_n(S_{g,p})$ by acting on coordinates. Let $\sigma(f): \text{PConf}_n(S_{g,p})\to \text{PConf}_n(S_{g,p})$ be the permutation of coordinates according to $T(f)$. We have that the permutation $\sigma(f)^{-1}\circ f_n$ fixes the base point $(b_1,...,b_n)$, which induces a map on the fundamental group $(\sigma(f)^{-1}\circ f_n)_*:PB_n(S_{g,p})\to PB_n(S_{g,p})$. Geometrically, we have the identification $C(f)=(\sigma(f)^{-1}\circ f_n)_*$. 
 
Let $F_i: PB_n(S_{g,p})\to  PB_{n-1}(S_{g,p})$ be the forgetful map forgetting the $i$th coordinate. By Theorem \ref{main1}, given any automorphism $A:PB_n(S_{g,p})\to PB_n(S_{g,p})$, we have that $F_i\circ A$ factors through $F_j$ for some $j\in \{1,...,n\}$. Therefore we obtain a homomorphism $S:\text{Aut}(PB_n(S_g))\to \Sigma_n$, where for $A\in \text{Aut}(PB_n(S_g))$, the image $S(A)$ satisfies that $F_i\circ A$ factors through $F_{S(A)(i)}$. Since $C(f)=(\sigma(f)^{-1}\circ f_n)_*$ and $\sigma(f)^{-1}\circ f_n$ permutes coordinates of $\text{PConf}_n(S_{g,p})$ by $T(f)^{-1}$, we know that $T(f)^{-1}=S(C(f))$.

By the five lemma, to prove that $C$ is an isomorphism, we only need to show that $C':\text{Ker}(T)\to \text{Ker}(S)$ is isomorphism. Denote by $\text{Aut}^0(PB_n(S_g))$ the subgroup of $\text{Aut}(PB_n(S_g))$ consisting of all $A\in \text{Aut}(PB_n(S_{g,p}))$ such that $F_n\circ A$ factors through $F_n$; i.e. $\text{Aut}^0(PB_n(S_g))=S^{-1}(\Sigma_{n-1}\times 1)$. We have that $T^{-1}(\Sigma_{n-1}\times 1)=\text{Mod}^\pm(S_{g,p,n-1,1})$. By the five lemma, to show that $C'$ is an isomorphism, we only need to show that $C'': S^{-1}(\Sigma_{n-1}\times 1)\to T^{-1}(\Sigma_{n-1}\times 1)$ is an isomorphism; i.e. $C'': \text{Mod}^\pm(S_{g,p,n-1,1})\to \text{Aut}^0(PB_n(S_g))$ is an isomorphism. The action $A\in \text{Aut}^0(PB_n(S_g))$ factors through $F_n$, which induces an action of the following commutative diagram:
\[
\xymatrix{
1\ar[r] &\pi_1(S_{g,p+n-1}) \ar[r]\ar[d]^R & PB_n(S_{g,p})\ar[r]^{F_n}\ar[d]^A & PB_{n-1}(S_{g,p})\ar[r]\ar[d]^{} & 1\\
1\ar[r] &\pi_1(S_{g,p+n-1}) \ar[r] & PB_n(S_{g,p})\ar[r]^{F_n} & PB_{n-1}(S_{g,p})\ar[r] & 1.}
\]
Using this diagram, we can define a homomorphism $P:\text{Aut}^0(PB_n(S_{g,p}))\to \text{Aut}(\pi_1(S_{g,p+n-1}))$ via the formula $P(A)=R$. We will prove that $C''$ is an isomorphism by showing that $P$ is the inverse of $C''$. This will conclude the proof of Theorem \ref{main10}. 
	\begin{itemize}
		\item Step 1: $P$ is injective. 
		
		Consider $A\in \text{Ker}(P)$. Therefore $R=P(A)$ is the identity. For any $e\in PB_n(S_{g,p})$ and $x\in \pi_1(S_{g,p+n-1})$, we have that		\[
		exe^{-1}\xlongequal{R=id}R(exe^{-1})\xlongequal{R=A}A(e)R(x)A(e)^{-1}\xlongequal{R=id}A(e)xA(e)^{-1}		\]
		So $e^{-1}A(e)\in PB_n(S_{g,p})$ commutes with $x\in \pi_1(S_{g,p+n-1})$. Since the conjugation action of $PB_n(S_{g,p})$ on $\pi_1(S_{g,p+n-1})$ is faithful, this implies that $e^{-1}A(e)=id$. So $A(e)=e$, as desired. 
		\item Step 2: Image$(P)< \phi(\text{Mod}^{\pm}(S_{g,p,n-1,1}))$. 
		
		Consider $A\in \text{Aut}^0(PB_n(S_{g,p}))$ and set $R=P(A)$. For $e\in PB_n(S_{g,p})$, denote by $\sigma_e$ the automorphism of $\pi_1(S_{g,p+n-1})$ induced by $e$ via conjugation. For any $e\in PB_n(S_{g,p})$ and $x\in \pi_1(S_{g,p+n-1})$, we have
		\[
		R(exe^{-1})=A(e)R(x)A(e)^{-1}.\]
		That is to say 
		\[
		R\sigma_eR^{-1}=\sigma_{A(e)} \in \text{Aut}\big(\pi_1(S_{g,p+n-1})\big).\]
		By Proposition \ref{haha}, we have that $R\in \text{Mod}^{\pm}(S_{g,p,n-1,1})$. 
\item Step 3: Both $\pi_1(S_{g,p+n-1})$ and $PB_n(S_{g,p})$ are normal subgroups of $\text{Mod}^{\pm}(S_{g,p+n-1,1})$. For $f\in \text{Mod}^{\pm}(S_{g,p,n-1,1})$, the image $\phi(f)$ is the conjugation action of $f$ on $\pi_1(S_{g,p+n-1})$. The image $A=C''(f)$ is the conjugation action of $f$ on $PB_n(S_{g,p})$. The image $P\circ C''(f)\in \text{Aut}^0(\pi_1(S_{g,p+n-1}))$ is the restriction of $C''(f)$ on $\pi_1(S_{g,p+n-1})$, which is equal to the direct conjugation $\phi(f)$. Therefore we have that $\phi=P\circ C''$. Therefore the composition $P\circ C'': \text{Mod}^{\pm}(S_{g,p,n-1,1})\to \text{Aut}^0(\pi_1(S_{g,p+n-1}))$ is an isomorphism onto the image $\phi(\text{Mod}^{\pm}(S_{g,p,n-1,1}))$. By step 2, this means that $P$ is surjective onto $\phi(\text{Mod}^{\pm}(S_{g,p,n-1,1}))$. By step 1, we know that $P$ is also injective. So $P$ is an isomorphism, as desired. This implies that $C''$ is an isomorphism as well.
	\end{itemize}
The result for $B_n(S)$ follows from a theorem of Ivanov \cite[Theorem 2]{ivanovpermutation} saying that $PB_n(S)$ is a characteristic subgroup of $B_n(S)$. 
\end{proof}

\bibliography{citing}{}

\end{document}